\documentclass[a4paper,12pt]{article}
\usepackage[cp1251]{inputenc}
\usepackage[english]{babel}
\usepackage{amsthm,amsmath,amssymb,latexsym}
\begin{document}
\begin{center}
\large{\textbf{Distributive Invariant Centrally Essential Rings}}
\end{center}

\hfill {\sf Askar Tuganbaev}

\hfill National Research University "MPEI"

\hfill e-mail: tuganbaev@gmail.com

\textbf{Abstract.} In recent years, centrally essential rings have been intensively studied in ring theory. In particular, they find applications in homological algebra, group rings, and the structural theory of rings. The class of essentially central rings strongly extends the class of commutative rings. For such rings, a number of recent papers contain positive answers to some important questions from ring theory that previously had positive answers for commutative rings and negative answers in the general case. This work is devoted to a similar topic.
A familiar description of right Noetherian, right distributive centrally essential rings is generalized on a larger class of rings. Let $A$ be a ring with prime radical $P(A)$. It is proved that $A$ is a right distributive, right invariant centrally essential ring and $P(A)$ is a finitely generated right ideal such that the factor-ring $A/P(A)$ does non contain an infinite direct sum of non-zero ideals if and only if $A=A_1\times\cdots\times A_n$, where every ring $A_k$ is either a commutative Pr\"ufer domain or an Artinian uniserial ring.


\section{Introduction}
We consider only associative rings with non-zero identity element. The phrases of the form '$A$ is a right (resp., left) Artinian ring' mean that the module $A_A$ (resp., $_AA$) is Artinian. The phrases of the form '$A$ is an Artinian (resp., uniserial) ring' mean that $A_A$ and $_AA$ are Artinian (resp., uniserial) modules.

A module $M$ is said to be \textsf{invariant} if all submodules of $M$ are fully invariant in $M$. It is clear that a ring $A$ is right (resp., left) invariant if and only if every right (resp., left) ideal of $A$ is an ideal.

A module $M$ is said to be \textsf{distributive} if the lattice of all submodules of $M$ is distributive. A module $M$ is said to be \textsf{uniserial} if any two submodules of $M$ are comparable with respect to inclusion. It is clear that every uniserial module is distributive. The ring of integers $\mathbb{Z}$ is an example of a commutative distributive non-uniserial ring.

A ring $A$ with center $C$ is said to be \textsf{centrally essential} if either the ring $A$ is commutative or for any non-zero non-central element $a\in A$ there exist non-zero central elements $x,y\in C$ with $ax=y$, i.e., the module $A_C$ is an essential extension of the module $A_C$.

\textbf{Remark 1.1.} It is clear that all commutative rings are centrally essential. There exist non-commutative, centrally essential, finite rings.\\ Let $F$ be a field consisting of three elements, $V$ be a linear $F$-space with basis $e_1,e_2,e_3$, and let $\Lambda(V)$ be the Grassmann algebra of the space $V$. Since $e_1\wedge e_1=e_2\wedge e_2=e_3\wedge e_3=0$ and any product of generators is equal to the $\pm$product of generators with ascending subscripts, $\Lambda(V)$ is a finite $F$-algebra of dimension 8 with basis
$$
\{1,e_1,e_2,e_3,e_1\wedge e_2,e_1\wedge e_3,e_2\wedge e_3,e_1\wedge e_2\wedge e_3\},
$$
$$
|\Lambda(V)|=3^8,\;e_k\wedge e_i\wedge e_j=-e_i\wedge e_k\wedge e_j=e_i\wedge e_j\wedge e_k.
$$
Therefore, if
$$
x=\alpha_0\cdot 1+\alpha_1^1e_1+\alpha_1^2e_2+\alpha_1^3e_3+\alpha_2^1e_1\wedge e_2+\alpha_2^2e_1\wedge e_3+$$
$$
+\alpha_2^3e_2\wedge e_3+\alpha_3e_1\wedge e_2\wedge e_3,
$$ 
then 
$$
\begin{array}{lll}[e_1,x]&=&2\alpha_1^2e_1\wedge e_2+2\alpha_1^3e_1\wedge e_3,\\
{}[e_2,x]&=-&2\alpha_1^1e_1\wedge e_2+2\alpha_1^3e_2\wedge e_3,\\
{}[e_3,x]&=-&2\alpha_1^1e_1\wedge e_3-2\alpha_1^2e_2\wedge e_3.\end{array}
$$
Thus, $x\in Z(\Lambda(V))$ if and only if $\alpha_1^1=\alpha_1^2=\alpha_1^3=0.$ In other words, the center of the algebra $\Lambda(V)$ is of dimension $5$. On the other hand, if $\alpha_1^1\neq 0$, then
$$
x\wedge(e_2\wedge e_3)=\alpha_0e_2\wedge e_3+\alpha_1^1e_1\wedge e_2\wedge e_3\in Z(\Lambda(V))\setminus \{0\}.
$$ 
In addition, $e_2\wedge e_3\in Z(\Lambda(V))$. A similar argument applies if $\alpha_1^2\neq 0$ or $\alpha_1^3\neq 0$. Consequently, $\Lambda(V)$ is a finite centrally essential non-commutative ring.

The following result is proved in \cite{MT23}.

\textbf{Theorem 1.2; \cite[Theorem 1.2]{MT23}.} A centrally essential ring $A$ is a right distributive, right Noetherian ring if and only if $A=A_1\times\cdots\times A_n$, where every ring $A_k$ is either a commutative Dedekind domain or a (not necessarily commutative) Artinian uniserial ring.

\textbf{Remark 1.3.} Let $A$ be a right Noetherian ring. It is well known that every right ideal is finitely generated and each factor ring of $A$ does non contain an infinite direct sum of non-zero right ideals. If the right Noetherial ring $A$ is right distributive, as well, then $A$ is right invariant, \cite{Ste74}. In addition, every commutative Dedekind domain is a Pr\"ufer domain. 

It follows from Remark 1.3 that Theorem 1.2 is a corollary of the following Theorem 1.4 which is the main result of this paper. 

\textbf{Theorem 1.4.} For a ring $A$ with prime radical $P$, the following conditions are equivalent.
\begin{enumerate}
\item[\textbf{1)}]
$A$ is a right distributive, right invariant, centrally essential ring and $P$ is a finitely generated right ideal such that the factor ring $A/P$ does non contain an infinite direct sum of non-zero ideals.
\item[\textbf{2)}]
$A=A_1\times\cdots\times A_n$, where every ring $A_k$ is either a commutative Pr\"ufer domain or an Artinian uniserial ring.
\end{enumerate}

\textbf{Remark 1.5.} In connection to Theorems 1.4 and 1.2, we note that there exist rings which satisfy Theorem 1.4 but do not satisfy Theorem 1.2. Let $F$ be a field, $G$ be the multiplicative group of positive rational numbers, and let $A$ be the group ring $FG$. It is easy to verify that $A$ is a commutative non-Noetherian domain in which every finitely generated ideal is principal. Therefore, $A$ is a Pr\"ufer domain which is not a Dedekind domain. Thus,
$A$ satisfies Theorem 1.4 but does not satisfy Theorem 1.2.

The proof of Theorem 1.4 is given in the next section; the proof is based on several assertions, some of which are of independent interest.

We give some definitions. Direct summands of free modules are called \textsf{projective} modules. A module $M$ is said to be \textsf{semihereditary} if all finitely generated submodules of the module $M$ are projective. A module $M$ is said to be \textsf{hereditary} if all submodules of the module $M$ are projective. 

A ring $A$ is called a \textsf{domain} if $A$ does not have non-zero zero-divisors. A commutative domain $A$ is called a \textsf{Pr\"ufer} domain if $A$ is a commutative semihereditary domain. A commutative domain $A$ is called a \textsf{Dedekind} domain if $A$ is a commutative hereditary domain. 

A module $M$ is said to be \textsf{Noetherian} (resp., \textsf{Artinian}) if $M$ does not contain an infinite properly ascending (resp., properly descending) chain of submodules. 
A module $M$ is said to be \textsf{finite-dimensional} (in the sense of Goldie) if $M$ does not contain a submodule which is an infinite direct sum of non-zero submodules. If $A$ is a ring, then a proper ideal $B$ of $A$ is said to be \textsf{completely prime} if the factor ring $A/B$ is a domain.

Other definitions, used in the paper, are standard; e.g., see \cite{Fai73}, \cite{Wis91}.

\section{Proof of Theorem 1.4}\label{section2}

\textbf{Lemma 2.1; see \cite[Corollary 9.13]{Fai73} and \cite[Theorem 2.2]{Ste75}.}\\
Let $A$ be a right finite-dimensional, right non-singular, semiprime ring. Then every essential right ideal of $A$ contains a non-zero-divisor.

\textbf{Lemma 2.2.} Let $A$ be a right distributive ring with prime radical $P$.

\textbf{a; \cite{Ste74}.} For any two elements $x,y\in A$, there exist elements $a,b\in A$ such that $a+b=1$ and $xaA+ybA\subseteq xA\cap yA$. In addition, if $xA\cap yA=0$, then $\text{Hom}_A(xA,yA)=\text{Hom}_A(yA,xA)=0$, $AyAxA=AxAyA=0$, and there exist elements $a,b\in A$ such that $a+b=1$, $xaA=ybA=0$, and $xA\oplus yA=(x+y)A$. 

\textbf{b; \cite{Ste74}.} In the right distributive rings $A$ and $A/P$, all idempotents are central. In addition, if the ring $A$ is right Noetherian, then $A$ is right invariant.

\textbf{c.} The ring $A$ is indecomposable (i.e., $A$ does non have non-trivial central idempotents) if and only if the ring $A/P$ is indecomposable, if and only if $A$ does non have non-trivial idempotents, if and only if $A/P$ does not have non-trivial idempotents. In addition, if the ring $A$ is local, then $A$ is a right uniserial ring; \cite{Ste74}.

\textbf{d.} If $A$ has a completely prime nil-ideal $M$, then $M=xM$ for any element $x\in A\setminus M$.

\textbf{e.} If $A/P$ is a right finite-dimensional, right non-singular,  indecomposable ring, then $P$ is a completely prime nil-ideal and $P=xP$ for any element $x\in A\setminus P$.

\begin{proof} The assertions \textbf{a} and \textbf{b} are partial cases of \cite[Theorem 1.6, Corollaries 1,2 of Proposition 1.1]{Ste74}. For convenience, we give brief proofs of \textbf{a} and \textbf{b}.

\textbf{a.} Let $T=xA\cap yA$. Since
$(x+y)A=(x+y)A\cap xA+(x+y)A\cap yA$, there exist $b,d\in A$ such that
$$
(x+y)b\in xA,\quad (x+y)d\in yA,\quad x+y=(x+y)b+(x+y)d.
$$
Therefore, $yb=(x+y)b-xb\in T$ and $xd=(x+y)d-yd\in T$. We set $a= 1-b$ and $z= a-d=1-b-d$. Then
$$
1=a+b,\; (x+y)z=(x+y)-(x+y)b-(x+y)d=0,
$$
$$
xa=xd+xz=xd+(x+y)z-yz=xd-yz,\;
yz=-xz\in T, \; xa\in T.
$$
Now we assume that $xA\cap yA=0$, it follows from the above that there exist elements $a,b\in A$ such that $a+b=1$ and $xaA=ybA=0$. Let  $c\in A$, $x'=xc\in xA$, $f\in\text{Hom}_A(xA,yA)$ and $y'=f(x')\in f(xA)$. By the above,  there exist elements $a',b'\in A$ such that $a'+b'=1$ and $x'a'A=y'b'A=0$. Then
$$
y'=f(x')=f(x'a')+f(x'b')=f(x'b')=y'b'=0.
$$
Therefore, $\text{Hom}_A(xA,yA)=0$. For every $c\in A$, we have a homomorphism 
$f\in \text{Hom}_A(xA,yA)$ defined by the relation $f(xa)=yca$. Since $\text{Hom}_A(xA,yA)=0$, we have $yca=0$. Therefore, $AyAxA=0$. Similarly, we have $\text{Hom}_A(yA,xA)=0$ and $AxAyA=0$. 

It remains to prove that $xA\oplus yA=(x+y)A$. Since $(x+y)A\subseteq xA\oplus yA$, it is sufficient to prove that $x,y\in (x+y)A$. Indeed, $x=xa+xb=xb=(x+y)b$ and $y=ya+yb=ya=(x+y)a$.

\textbf{b.} Since $eA\cap (1-e)A=0$, it follows from \textbf{a} that $(1-e)Ae=eA(1-e)=0$. Therefore, the idempotent $e$ is central. The second assertion is proved in \cite{Ste74}.

\textbf{c.} Since all idempotents are lifted modulo the nil-ideal $P$, the first assertion follows from \textbf{b}.

Let the ring $A$ be local, $x,y\in A$, and let $x\notin yA$. By \textbf{a}, there exist elements $a,b\in A$ such that $a+b=1$, $xaA\subseteq yA$, and $ybA\subseteq xA$. Since $x\notin yA$, the element $a$ is not invertible. Therefore, the element $b$ is invertible, since $a+b=1$ and the ring $A$ is local. Then $y=(yb)b^{-1}\in xA$ and the ring $A$ is right uniserial.

\textbf{d.} Let $x\in A\setminus M$ and let $y$ be an arbitrary element of the ideal $M$. By \textbf{a}, there exist two elements $a,b\in A$ such that $a+b=1$, $xa\in yA$ and $yb\in xA$. The ideal $M$ is completely prime, $x\in A\setminus M$ and $xa\in M$. Therefore, $a\in M$ and the element $a$ is nilpotent. Therefore, the element $b=1-a$ is invertible; in addition, $yb\in xA$. Then $y=xz$ for some $z\in A$. Since the element $xz$ is contained in the completely prime ideal $M$ and $x\in A\setminus M$, we have that $z\in M$ and $y=xz\in xM$. Since the element $y\in M$ is arbitrary, we have $M=xM$.

\textbf{e.} We set $R=A/P$. Since $P$ is a nil-ideal, it follows from \textbf{d} that it is sufficient to prove that $R$ is a domain. Since the ring $A$ is indecomposable, it follows from \textbf{b} that the right distributive, right finite-dimensional, right non-singular, semiprime ring $R$ does not have non-trivial idempotents. Therefore, the right $R$-module $R$ is indecomposable. 
Since $R$ is right finite-dimensional, the ring $R$ has an essential right ideal $B$ of the form $B_1\oplus\ldots\oplus B_n$ for some non-zero uniform right ideals $B_k$. By Lemma 2.1, the essential right ideal $B$ contains some non-zero-divisor $b=b_1+\ldots + b_n$, where $0\ne b_k\in B_k$ for some $k$. Then $R_R\cong bR$, whence the indecomposable module $R_R$ is a direct sum of uniform modules. Therefore, the right non-singular ring $R$ is right uniform, whence $R$ is a domain, which is required.\end{proof} 

\textbf{Lemma 2.3.} For a ring $A$ with prime radical $P$, the following conditions are equivalent.
\begin{enumerate}
\item[\textbf{1)}]
$A$ is a right distributive ring and $A/P$ is a right finite-dimensional, right non-singular ring.
\item[\textbf{2)}]
$A=A_1\times\cdots\times A_n$, where every ring $A_k$ is a right distributive ring with prime radical $P_k$ such that $P_k$ is a completely prime nil-ideal of $A_k$ and $P_k=xP_k$ for any element $x\in A_k\setminus P_k$.
\end{enumerate}
\begin{proof} 1)\,$\Rightarrow$\,2). The right finite-dimensional, right non-singular  ring $A/P$ is the finite direct product of indecomposable, right finite-dimensional, right non-singular rings. Without loss of generality, we can assume that the ring $A$ is indecomposable. By Lemma 2.2(e), we have that 
$P$ is a completely prime nil-ideal and $P=xP$ for any element $x\in A\setminus P$.

2)\,$\Rightarrow$\,1). The assertion is directly verified.\end{proof} 

\textbf{Lemma 2.4.} For a ring $A$ with prime radical $P$, the following conditions are equivalent.
\begin{enumerate}
\item[\textbf{1)}]
$A$ is a right distributive, right invariant ring and the factor ring $A/P$  does non contain an infinite direct sum of non-zero ideals.
\item[\textbf{2)}]
$A=A_1\times\cdots\times A_n$, where every ring $A_k$ is a right distributive, right invariant ring with prime radical $P_k$ such that $P_k$ is a completely prime nil-ideal of $A_k$ and $P_k=xP_k$ for any element $x\in A_k\setminus P_k$.
\end{enumerate}
\begin{proof}  1)\,$\Rightarrow$\,2). We set $R=A/P$. Since the right invariant ring $R$ does non contain an infinite direct sum of non-zero ideals, $R$ is right finite-dimensional. By Lemma 2.3, it is sufficient to prove that the ring $R$ is right non-singular ring. Let $a$ be an element of $R$ whose right annihilator $B$ is an essential right ideal. We have to prove that $a=0$. Assume the contrary. Then the square of the non-zero ideal $B\cap aR$ of the right invariant semiprime ring is equal to zero; this is a contradiction.

2)\,$\Rightarrow$\,1). The assertion is directly verified with the use of Lemma 2.3.\end{proof} 

\textbf{Lemma 2.5.; see \cite[Lemma 20]{Tug95}.} If $A$ is a right invariant ring and $M$ is a finitely generated right $A$-module, then for any submodule $N$ of $M$, there exists an ideal $D$ of $A$ with $MD=N$.

\textbf{Lemma 2.6; \cite[Corollary 3]{Jen66}.} A commutative domain is distributive if and only if it is a Pr\"ufer domain.

\textbf{Proposition 2.7.} Let $A$ be a right distributive, right invariant,  indecomposable ring with prime radical $P$, the right ideal $P$ be finitely generated, and let the factor ring $A/P$ do not contain an infinite direct sum of non-zero ideals. 
\begin{enumerate}
\item[\textbf{a.}]
$A/P$ is a right distributive right invariant domain, $P$ is a completely prime nilpotent ideal and $P=xP$ for any element $x\in A\setminus P$.
\item[\textbf{b.}]
For each submodule $N$ of $P_A$, there exists an ideal $D$ of $A$ such that $N=MD$.
\item[\textbf{c.}]
For each submodule $N$ of $P_A$, there exists an ideal $D$ of $A$ such that $N=xN=XN$ for any element $x\in A\setminus M$ and each ideal $X$ properly containing $P$.
\item[\textbf{d.}] 
If $P$ contains a non-zero central element $m$, then $A$ is a right uniserial right Artinian ring with radical $P$.
\item[\textbf{e.}] 
If the ring $A$ is centrally essential, then either $A$ is a commutative Pr\"ufer domain or $A$ is a uniserial Artinian ring.
\end{enumerate}

\begin{proof} \textbf{a.} By Lemma 2.4, the prime radical $P$ is a completely prime nil-ideal of $A$ and $P=xP$ for any element $x\in A_k\setminus P_k$. In particular, $A/P$ is a right distributive right invariant domain. Since the ring $A$ is right invariant and $P$ is a finitely generated right nil-ideal, the ideal $P$ is nilpotent.

\textbf{b.} The assertion follows from Lemma 2.5.

\textbf{c.} By \textbf{b}, there exists an ideal $D$ of $A$ with $PD=N$. By \textbf{a}, $P=xP$. Therefore, $N=(xP)D=x(PD)=xN$. In particular, $N=XN$ for  each ideal $X$ properly containing $P$.

\textbf{d.} Since the ring $A$ is right invariant and $m$ is a non-zero central element, there exists a maximal ideal $X$ of the ring $A$ such that $A/X$ is a division ring and the ideal $mA$ properly contains the ideal $(mA)X=X(mA)$. If the maximal ideal $X$ properly contains $P$, the it follows from \textbf{c} that $mA=X(mA)=(mA)X$; this is a contradiction. Therefore, $X=P$, whence $A$ is a local ring with nilpotent Jacobson radical $P$ and $A$ is a right uniserial ring, by Lemma 2.2(c). Now it is directly verified that $A$ is a right uniserial, right Artinian ring.

\textbf{e.} Assume that the completely prime ideal $P$ is equal to the zero. Then $A$ is a centrally essential domain. Let $a$ be any non-zero element of $A$. Then $ax=y$ for some non-zero central elements of $A$. Therefore, it is directly verified that $A$ is an \"Ore domain. Thus, $A$ is an order in some division ring $Q$, the central element $x$ of $A$ is invertible in $Q$, and $a=yx^{-1}$ is a central element of $A$. Therefore, the domain $A$ is commutative. By Lemma 2.6, the commutative distributive domain is a Pr\"ufer domain.

Assume that $P\ne 0$. Since the ring $A$ is centrally essential, the ideal $P$ contains a non-zero central element. It follows from \textbf{d} that $A$ is a right uniserial, right Artinian ring. By Theorem 1.2, the ring $A$ is a uniserial Artinian ring.\end{proof}

\textbf{The completion of the proof of Theorem 1.4.}\\
Let $A$ be a centrally essential ring.

If $A=A_1\times\cdots\times A_n$, where every ring $A_k$ is a either commutative Dedekind domain or (not necessarily commutative) Artinian uniserial ring, then the assertion follows from Lemma 2.4 and Lemma 2.6.

Now let $A$ be a right distributive, right invariant, centrally essential ring and $P$ be a finitely generated right ideal such that the factor ring $A/P$ does non contain an infinite direct sum of non-zero ideals.
Without loss of generality, we can assume that the ring $A$ is indecomposable. By Proposition 2.7(e), either $A$ is a commutative Pr\"ufer domain or $A$ is a uniserial Artinian ring.

Askar Tuganbaev is supported by Russian Scientific Foundation.

\end{document}